\documentclass[twoside]{article}
\newtheorem{theo}{Theorem}

\newtheorem{prop}[theo]{Proposition}
\newtheorem{rema}[theo]{Remark}
\def\al{\alpha}

\def\De{\Delta}
\def\beq{\begin{equation}}
\def\eeq{\end{equation}}
\def\bea{\begin{eqnarray}}
\def\eea{\end{eqnarray}}
\def\beas{\begin{eqnarray*}}
\def\eeas{\end{eqnarray*}}
\def\nn{\nonumber}
\def\ds{\displaystyle}
\def\U{{\cal U}_h({\rm sl}(2))}
\def\D{{\cal D}}
\begin{document}
\renewcommand{\thefootnote}{\fnsymbol{footnote}}
\setcounter{footnote}{1}

\begin{center}
{\bf THE JORDANIAN DEFORMATION OF SU(2) AND \\ 
CLEBSCH-GORDAN COEFFICIENTS}\footnote{Presented at the 6th
International Colloquium ``Quantum Groups and Integrable Systems,''
Prague, 19-21 June 1997.}    

\vspace{1cm}

{\sc Joris Van der Jeugt}\footnote{Research Associate of the
Fund for Scientific Research -- Flanders (Belgium); \\
E-mail~: Joris.VanderJeugt@rug.ac.be}\\[2mm]
{\sl Department of Applied Mathematics and Computer Science,} \\ 
{\sl University of Ghent, Krijgslaan 281-S9, B-9000 Gent, Belgium}
\end{center}

\vspace{1cm}

{\small Representation theory for the Jordanian quantum algebra $\U$ is
developed using a nonlinear relation between its generators and those
of ${\rm sl}(2)$. Closed form expressions are given for the action of
the generators of $\U$ on the basis vectors of finite dimensional
irreducible representations. In the tensor product of two such
representations, a new basis is constructed on which the generators of
$\U$ have a simple action. Using this basis, a general formula is
obtained for the Clebsch-Gordan coefficients of $\U$. Some remarkable
properties of these Clebsch-Gordan coefficients are derived.
}

\vspace{3mm}
\begin{center}
{\bf 1\ \ Introduction}
\end{center}
\nopagebreak
The group ${\rm GL}(2)$ admits, upto isomorphism, only two quantum
group deformations with central determinant~: ${\rm GL}_q(2)$ and ${\rm
GL}_h(2)$, see~[1]. The quantum group ${\rm GL}_q(2)$ has been well
studied, being the prototype example for many works on quantum groups.
Investigations of the Jordanian quantum group ${\rm GL}_h(2)$, or ${\rm
SL}_h(2)$, and its dual quantum algebra $\U$ started more recently.
Its defining relations were given in~[2,3], and a construction of the
dual Hopf algebra in~[4]. Recently, also for the 2-parameter Jordanian
quantum group  ${\rm GL}_{g,h}(2)$ its dual was constructed~[5].
For a development of its differential calculus or
differential geometry we refer to~[6] and~[7].
A construction of the universal $R$-matrix was given in~[8,9,10].

In this paper we are primarily interested in the irreducible 
finite dimensional representations of $\U$. Also here, there has been
progress in recent years. In~[11], a direct construction of these
representations was given by factorising the Verma module.
An important development was
given by Abdesselam {\em et al}~[12]~: they gave a
nonlinear relation between the generators of $\U$ and the classical
generators of $sl(2)$. As a consequence they obtained
expressions for the action of the generators of $\U$ on basis vectors
of the finite dimensional irreducible representations. These
expressions were not always in closed form, and this was solved in~[13].
In~[14], finite and infinite dimensional representations of $\U$ are
constructed, and for the first time the tensor product of two
representations is considered, yielding some examples of Clebsch-Gordan
coefficients. The problem of determining Clebsch-Gordan coefficients
was then completely solved in~[13].

In the present paper we shall discuss a number of interesting
properties of the Clebsch-Gordan coefficients of $\U$, after recalling
some of the main results of~[13].

\vspace{3mm}
\begin{center}
{\bf 2\ \ ${\rm SL}_h(2)$ and $\U$}
\end{center}
\nopagebreak
Consider the bialgebra ${\cal A}_h(2)$ with parameter $h$ and four
generators $a$, $b$, $c$, $d$ subject to the relations~:
\beq
\begin{array}{ll}
ba=ab-ha^2+h\D &
ca=ac+hc^2\\
da=ad+hdc-hac &
bd=db-hd^2+h\D\\
cd=dc+hc^2 &
cb=bc+hdc+h^2c^2
\end{array}
\eeq
where $\D=ad-bc-hac$. It is easy to verify the $\D$ is central.
With $t=\left( {a \atop c}\ {b \atop d}\right)$, there is a
comultiplication given by $\Delta(t)=t \otimes t$, and a co-unit
$\epsilon(t) = \left( {1 \atop 0}\ {0 \atop 1}\right)$, turning ${\cal
A}_h(2)$ into a coalgebra. The element
$\D$ is group-like, so one can extend ${\cal A}_h(2)$ by $\D ^{-1}$,
and then an antipode $S$ can be defined leading to the Hopf algebra
${\rm GL}_h(2)$.
Putting $\D = 1$ gives rise to the matrix quantum group ${\rm
SL}_h(2)$, see~[11].

The dual Hopf algebra of ${\rm SL}_h(2)$ is denoted by $\U$. It is
an associative algebra generated by $H$, $Y$, $T$ and $T^{-1}$
satisfying quadratic relations~[4]. For us it is more convenient to work with
$X=({\rm log}T)/h$, i.e.\ $T=e^{hX}$ and $T^{-1}=e^{-hX}$. Then the
relations read~:
\bea
&&[H,X]=2 {\sinh hX \over h}, \qquad \qquad [X,Y]=H, \nn\\
&&[H,Y]= -Y(\cosh hX) - (\cosh hX)Y. 
\eea
The comultiplication is given by~:
\bea
&&\De(H)=H\otimes e^{hX}+ e^{-hX}\otimes H, \nn\\
&&\De(X)=X\otimes 1 + 1\otimes X, \\
&&\De(Y)=Y\otimes e^{hX}+ e^{-hX}\otimes Y. \nn
\eea
The other ingredients (co-unit, antipode) are also defined, but not
needed here.

\vspace{3mm}
\begin{center}
{\bf 3\ \ Relation between $\U$ and ${\rm sl}(2)$, and representations}
\end{center}
\nopagebreak
With the following definition~[12]
\bea
&& Z_+ = {2 \over h} \tanh {hX \over 2},\nn\\
&& Z_- = ( \cosh {hX \over 2}) Y (\cosh {hX\over 2}), 
\eea
it follows that the elements $\{H, Z_+, Z_-\}$ satisfy the commutation
relations of a classical ${\rm sl}(2)$ basis~:
\[
[H,Z_\pm]=\pm 2 Z_{\pm}, \qquad [Z_+,Z_-]=H.
\]
These relations can be inverted, e.g.\ 
\[
e^{hX} = (1+{h\over 2} Z_+)   (1-{h\over 2} Z_+)^{-1}.
\]
These relations can also be used to give explicit matrix elements for the
finite dimensional representations of $\U$.

Recall that finite dimensional irreducible representations of ${\rm sl}(2)$
are labeled by a number $j$, with
$2j$ a non-negative integer. The representation space can be denoted by
$V^{(j)}$ with basis $e^j_m$ ($m= -j, -j+1,\ldots, j$), and the action is
\bea
&& H e^j_m = 2m\; e^j_m, \nn\\
&& Z_\pm e^j_m = \sqrt{ (j\mp m)(j\pm m+1) }\; e^j_{m\pm 1}.
\eea
For us, a more convenient basis for computations is the following
$v$-basis related to the above $e$-basis by~:
\[
v^j_m = \al_{j,m} e^j_m,\;\;\hbox{with } \al_{j,m} = \sqrt{
(j+m)!/(j-m)!}\; .
\]
The ${\rm sl}(2)$ matrix elements in this basis are~:
\bea
&& H v^j_m = 2m\; v^j_m, \nn\\
&& Z_+ v^j_m = v^j_{m+1},  \\
&& Z_- v^j_m = (j+ m)(j- m+1)\; v^j_{m-1}., \nn
\eea
where $v^j_{j+1}\equiv 0$.

Using the explicit mapping between $\{H,Z_+,Z_-\}$ and $\{H,X,Y\}$,
plus a number of combinatorial identities~[13], we obtained~:

\begin{prop}
The action of the generators of $\U$ on
the representation space $V^{(j)}$ is given by
\bea
&& H v^j_m = 2m\; v^j_m, \nn\\
&& X v^j_m = \sum_{k=0}^{\lfloor (j-m-1)/2 \rfloor} {(h/2)^{2k} \over
2k+1} \; v^j_{m+1+2k} , \\
&&Y v^j_m = (j+m)(j-m+1) v^j_{m-1} - (j-m)(j+m+1)\left(h\over 2\right)^2
v^j_{m+1}  \nn\\
&&\qquad +\sum_{s=1}^{\lfloor(j-m+1)/2\rfloor} \left(h\over
2\right)^{2s} v^j_{m-1+2s}, \nn
\eea
\end{prop}
It should be noted that the matrix elements of $X$ were already
obtained in~[12]. Those of $Y$ were also determined in~[12], however not in
closed form but as a complicated sum. In~[13] we showed how such sums
can be reduced to a simple form, using recently developed
algorithms~[15]. Proposition~1 is easy to apply and gives immediately
all matrix elements of the $\U$ generators. For example, the
representatives for $X$ and $Y$, respectively, in the $v$-basis
for $j=2$ are given by~:
\[
\left(\begin{array}{ccccc}
0 & 1 & 0 & h^2/12 & 0 \\
0 & 0 & 1 & 0 & h^2/12 \\
0 & 0 & 0 & 1 & 0 \\
0 & 0 & 0 & 0 & 1 \\
0 & 0 & 0 & 0 & 0 
\end{array}\right) \quad
\left(\begin{array}{ccccc}
0 & -3h^2/4 & 0 & h^4/16 & 0 \\
4 & 0 & -5h^2/4 & 0 & h^4/16 \\
0 & 6 & 0 & -5h^2/4 & 0 \\
0 & 0 & 6 & 0 & -3h^2/4 \\
0 & 0 & 0 & 4 & 0 
\end{array}\right);
\]
with $H$ given by the usual matrix ${\rm diag}(4,2,0,-2,-4)$.
Note that the ${\rm sl}(2)$ representatives in the $v$-basis are
recovered simply by putting $h=0$.

\vspace{3mm}
\begin{center}
{\bf 4\ \ Tensor product of $\U$ representations}
\end{center}
\nopagebreak
Consider $V^{(j_1)}\otimes V^{(j_2)}$ with basis $v^{j_1}_{m_1} \otimes
v^{j_2}_{m_2}$. Our purpose is to show that this decomposes into the direct
sum of representations $V^{(j)}$, $j=|j_1-j_2|,\ldots, j_1+j_2$. Note
that the vectors $v^{j_1}_{m_1} \otimes v^{j_2}_{m_2}$ are in general
no eigenvectors of $\De(H)$, since the comultiplication is given by~:
\bea
\De(H) &=& H\otimes e^{hX} + e^{-hX} \otimes H \nn\\
&=& H\otimes 1 + 1 \otimes H + 2H \otimes \sum_{n=1}^\infty
\left(hZ_+\over 2 \right)^n +  \sum_{n=1}^\infty
\left(- hZ_+\over 2 \right)^n \otimes 2H.
\eea
The eigenvectors of $\De(H)$ are linear combinations of the vectors
$v^{j_1}_{m_1} \otimes v^{j_2}_{m_2}$, and the coefficients play a
crucial role in this work. To define these coefficients, recall the
definition of the Pochhammer symbol~:
\beq
(a)_n = \left\{ 
 \begin{array}{lll} 
 a(a+1)\cdots (a+n-1) & \hbox{if} & n=1,2,\ldots; \\
 1 & \hbox{if} & n=0.
 \end{array} \right.
\eeq
Next we define
\beq
b^{m_1,m_2}_{k,l} = 
\left\{ 
 \begin{array}{ll}
\ds{(-2m_1-k)_l (-2m_2-l)_k \over k!l!} & \hbox{if }k\geq 0\hbox{ and }
l\geq 0 ; \\[2mm]
0 & \hbox{otherwise},
 \end{array} \right.
\label{cob}
\eeq
and finally the essential coefficients~:
\beq
a^{m_1,m_2}_{k,l} = (-1)^k (h/2)^{k+l} (b^{m_1,m_2}_{k,l} -
 b^{m_1,m_2}_{k-1,l-1} ).
\label{coa}
\eeq
Then we have the following important result~:

\begin{prop}
In $V^{(j_1)}\otimes V^{(j_2)}$, the vectors
\beq
w^{j_1,j_2}_{m_1,m_2} = \sum_{k=0}^{j_1-m_1} \sum_{l=0}^{j_2-m_2}
a^{m_1,m_2}_{k,l} \; v^{j_1}_{m_1+k} \otimes v^{j_2}_{m_2+l}
\label{wvector}
\eeq
form a basis consisting of eigenvectors of $\De(H)$. The explicit
action of $\De(H)$, $\De(X)$ and $\De(Y)$ is given by
\bea
&&\De(H)\, w^{j_1,j_2}_{m_1,m_2} = 2(m_1+m_2) \; w^{j_1,j_2}_{m_1,m_2},\nn\\
&&\De(Z_+)\, w^{j_1,j_2}_{m_1,m_2} = w^{j_1,j_2}_{m_1+1,m_2} +
w^{j_1,j_2}_{m_1,m_2+1}, \\
&&\De(Z_-)\, w^{j_1,j_2}_{m_1,m_2} = (j_1+m_1)(j_1-m_1+1)\; 
w^{j_1,j_2}_{m_1-1,m_2} + \nn\\
&& \qquad\qquad (j_2+m_2)(j_2-m_2+1)\; w^{j_1,j_2}_{m_1,m_2-1}.\nn
\eea
\end{prop}

\begin{rema}
{\rm This proposition tells us that the action of $\De(H)$,
$\De(X)$ and $\De(Y)$ on the $w$-vectors is the same as the action of
the ${\rm su}(2)$ generators (under the trivial Lie algebra
comultiplication) on the uncoupled vectors $v^{j_1}_{m_1} \otimes
v^{j_2}_{m_2}$. This observation implies the results on the tensor
product decomposition and Clebsch-Gordan coefficients for $\U$. In
particular, the Clebsch-Gordan coefficients for $\U$ are essentially
given by linear combinations of ${\rm su}(2)$ Clebsch-Gordan
coefficients, with $a^{m_1,m_2}_{k,l}$ the coefficients of this linear
combination.} 
\end{rema}

Let us first consider an example, say $V^{(1)}
\otimes V^{(1/2)}$. Using the formulas~(\ref{cob})-(\ref{wvector}), the
$w$-vectors are explicitly given by 
\[
\left( \begin{array}{c} w^{1,1/2}_{-1,-1/2} \\ w^{1,1/2}_{-1,1/2} \\
w^{1,1/2}_{0,-1/2} \\ w^{1,1/2}_{0,1/2} \\
w^{1,1/2}_{1,-1/2} \\ w^{1,1/2}_{1,1/2} \end{array}\right) =
\left( \begin{array}{cccccc}
1 & h & -h/2 & h^2/4 & h^2/4 & -h^3/8 \\[1.5mm]
0 & 1 & 0 & h/2 & 0 & 0 \\[1.5mm]
0 & 0 & 1 & 0 & -h/2 & h^2/4 \\[1.5mm]
0 & 0 & 0 & 1 & 0 & h/2 \\[1.5mm]
0 & 0 & 0 & 0 & 1 & -h \\[1.5mm]
0 & 0 & 0 & 0 & 0 & 1
\end{array} \right)
\left( \begin{array}{c} 
v^1_{-1}\otimes v^{1/2}_{-1/2} \\ v^1_{-1}\otimes v^{1/2}_{1/2} \\
v^1_{0}\otimes v^{1/2}_{-1/2} \\ v^1_{0}\otimes v^{1/2}_{1/2} \\
v^1_{1}\otimes v^{1/2}_{-1/2} \\ v^1_{1}\otimes v^{1/2}_{1/2} 
\end{array}\right).
\]
It is easy to verify that the inverse of the above upper-triangular
matrix is given by reflecting the matrix along its second diagonal,
i.e.\ by its skew-transpose~:
\[
\left( \begin{array}{cccccc}
1 & -h & h/2 & h^2/4 & 0 & -h^3/8 \\
0 & 1 & 0 & -h/2 & 0 & h^2/4 \\
0 & 0 & 1 & 0 & h/2 & h^2/4 \\
0 & 0 & 0 & 1 & 0 & -h/2 \\
0 & 0 & 0 & 0 & 1 & h \\
0 & 0 & 0 & 0 & 0 & 1
\end{array} \right).
\]
This turns out to be a general property of these matrices of
$a^{m_1,m_2}_{k,l}$ coefficients. In other words, we have

\begin{prop}
The coefficients $a^{m_1,m_2}_{k,l}$ satisfy
\beq
\sum_{n_1, n_2} a^{m_1,m_2}_{n_1-m_1,n_2-m_2}
a^{-M_1,-M_2}_{M_1-n_1,M_2-n_2} = \delta_{m_1,M_1} \delta_{m_2,M_2}.
\eeq
\end{prop}
Note that the above formula is nontrivial only for $M_1\geq m_1$ and
$M_2 \geq m_2$, otherwise the indices of the $a$-coefficients are
negative and thus automatically zero. The above property follows from
the following remarkable identity holding for arbitrary parameters $x$
and $y$~:
\bea
&&\sum_{k=0}^K \sum_{l=0}^L {(-x-k)_l(-y-l)_k \over k!l!}
{(x+K+k)_{L-l} (y+L+l)_{K-k} \over (K-k)!(L-l)!} {(xy+lx+ky)\over
(x+k)(y+l)} \nn\\
&&\times {(xy+Lx+Ky+lx+ky+2Kl+2kL)\over(x+K+k)(y+L+l)} = \delta_{K,0}
\delta_{L,0},
\label{xyident}
\eea
by putting $x=2m_1$, $y=2m_2$, $K=M_1-m_1$ and $L=M_2-m_2$. The proof of
(\ref{xyident}) falls beyond the scope of the present paper.

\vspace{3mm}
\begin{center}
{\bf 5\ \ Clebsch-Gordan coefficients and properties}
\end{center}
\nopagebreak
{}From Remark~3 it is easy to deduce that the decomposition of the
tensor product is given by
\[
V^{(j_1)} \otimes V^{(j_2)} = \bigoplus_{j=|j_1-j_2|}^{j_1+j_2}
V^{(j)},
\]
and we have

\begin{prop}
The Clebsch-Gordan coefficients for $\U$, in
\[
e^{(j_1j_2)j}_m = \sum_{n_1,n_2} {\cal C}^{j_1,j_2,j}_{n_1,n_2,m}(h) \;
e^{j_1}_{n_1} \otimes e^{j_2}_{n_2},
\]
are given by
\[
{\cal C}^{j_1,j_2,j}_{n_1,n_2,m}(h) = \sum_{m_1+m_2=m}
C^{j_1,j_2,j}_{m_1,m_2,m} A^{m_1,m_2}_{n_1-m_1,n_2-m_2},
\]
with $C^{j_1,j_2,j}_{m_1,m_2,m}$ the usual ${\rm su}(2)$ Clebsch-Gordan
coefficients, and $A^{m_1,m_2}_{n_1-m_1,n_2-m_2}$ determined by
\[
A^{m_1,m_2}_{k,l} = a^{m_1,m_2}_{k,l} {\al_{j_1,m_1+k}\al_{j_2,m_2+l} \over
\al_{j_1,m_1}\al_{j_2,m_2}}.
\]
\end{prop}

So apart from the $\al$-factors (which appear here because we have
formulated the proposition in the $e$-basis rather than in the
$v$-basis), the Clebsch-Gordan matrix is essentially the product of the
corresponding ${\rm su}(2)$ Clebsch-Gordan matrix with the upper
triangular matrix of $a$-coefficients considered in the previous
section. 

{}From the explicit form of the $a$-coefficients, and Proposition~5, it
follows that
\begin{prop}
The Clebsch-Gordan coefficients of $\U$ satisfy
\begin{itemize}
\item if $m=n_1+n_2$ then ${\cal
C}^{j_1,j_2,j}_{n_1,n_2,m}(h) = C^{j_1,j_2,j}_{n_1,n_2,m}$ ;
\item if $m>n_1+n_2$ then ${\cal C}^{j_1,j_2,j}_{n_1,n_2,m}(h) =0$;
\item if $m<n_1+n_2$ then ${\cal C}^{j_1,j_2,j}_{n_1,n_2,m}(h)$ is a
monomial in $h^{n_1+n_2-m}$.
\end{itemize}
\end{prop}

The most interesting property follows from Proposition~4~:
\begin{prop}
The Clebsch-Gordan coefficients of $\U$ satisfy the
skew-or\-tho\-go\-na\-lity relations
\beas
&&\sum_{n_1,n_2} (-1)^{j_1+j_2-j}{\cal C}^{j_1,j_2,j}_{n_1,n_2,m}(h) 
{\cal C}^{j_1,j_2,j'}_{-n_1,-n_2,-m'}(h) = \delta_{j,j'} \delta_{m,m'},\\
&&\sum_{j,m} (-1)^{j_1+j_2-j} {\cal C}^{j_1,j_2,j}_{n_1,n_2,m}(h) 
{\cal C}^{j_1,j_2,j}_{-n_1',-n_2',-m}(h) = \delta_{n_1,n_1'}
\delta_{n_2,n_2'} .
\eeas
\end{prop}

This property gives in fact the inverse matrix of a general
Clebsch-Gordan matrix of $\U$. The proof is as follows~: recall that
the Clebsch-Gordan matrix of $\U$ is essentially the product of an
upper-triangular matrix of $a$-coefficients with an ${\rm su}(2)$
Clebsch-Gordan matrix. But the upper-triangular matrix has an easy
inverse, namely its skew-transpose; and also the ${\rm su}(2)$
Clebsch-Gordan matrix has an easy inverse, namely its transpose (since
it is orthogonal). This, and some symmetry properties of ${\rm su}(2)$
Clebsch-Gordan coefficients, leads to Proposition~7.

\vskip 5mm
\begin{center}
{\bf References}
\end{center}
\nopagebreak
\begin{enumerate}

\item 
Kupershmidt B.\ A.: J.\ Phys.\ A~: Math.\ Gen. {\it 25} 
(1992), L1239.

\item 
Deminov E.\ E., Manin Yu.\ I., Mukhin E.\ E.\ and Zhdanovich D.~V.,
Progr.\ Theor.\ Phys.\ Suppl.\ {\it 102} (1990) 203.

\item 
Zakrzewski S.: Lett.\ Math.\ Phys. {\it 22} (1991), 287.

\item 
Ohn Ch.: Lett.\ Math.\ Phys. {\it 25} (1992), 85.

\item 
Aneva B.\ L., Dobrev V.\ K.\ and Mihov S.\ G.: Duality for the
Jordanian Matrix Quantum Group ${\rm GL}_{g,h}(2)$, preprint
q-alg/9705028. 

\item 
Aghamohammadi A.: Mod.\ Phys.\ Lett. {\it A8} (1993) 2607.

\item 
Karimipour V.: Lett.\ Math.\ Phys. {\it 30} (1994) 87.

\item 
Vladimirov A.\ A.: Mod.\ Phys.\ Lett. {\it A8} (1993) 2573.

\item 
Khorrami M., Shariati A., Abolhassani M.\ R.\ and Aghamohammadi~A.:
Mod.\ Phys.\ Lett. {\it A10} (1995) 873.

\item 
Ballesteros A.\ and Herranz F.\ J.: J.\ Phys.\ A~: Math.\ Gen. {\it 29}
(1996) L311.

\item 
Dobrev V.\ K.: Representations of the Jordanian quantum algebra $\U$,
preprint IC/96/14, ICTP, Trieste.

\item 
Abdesselam B., Chakrabarti A.\ and Chakrabarti R.: 
Mod.\ Phys.\ Lett.\ A {\it 11} (1996), 2883.

\item 
Van der Jeugt J.: Representations and Clebsch-Gordan coefficients for
the Jordanian quantum algebra $\U$, preprint q-alg/9703011.

\item 
Aizawa N.: Irreducible decompositon for tensor product representations
of Jordanian quantum algebras, preprint q-alg/9701022.

\item 
Petkov\v sek, M., Wilf, H.\ S.\ and Zeilberger, D.: {\em A=B}, A.\ K.\
Peters, Wellesley, Massachusetts, 1996.

\end{enumerate}
\end{document}